\documentclass[reqno,10pt]{amsart}

\usepackage{amssymb,graphicx,color}
\usepackage[latin1]{inputenc}

\theoremstyle{plain}
\newtheorem{thm}{Theorem}[section]

\newtheorem{prop}{Proposition}[section]

\newtheorem{defs}{Definition}[section]
\theoremstyle{definition}
\newtheorem{rmk}{Remark}[section]

\newcommand{\NN}{{\mathbb{N}}}

\newcommand{\RR}{{\mathbb{R}}}

\newcommand{\bu}{\mathbf{u}}
\newcommand{\bv}{\mathbf{v}}

\newcommand{\bx}{\mathbf{x}}

\newcommand{\bbf}{\mathbf{f}}

\newcommand{\bxi}{\boldsymbol{\xi}}
\newcommand{\bnabla}{\boldsymbol{\nabla}}

\newcommand{\define}{\stackrel{\text{\rm def}}{=}}
\newcommand{\Ccal}{{\mathcal C}}

\newcommand{\Scal}{{\mathcal S}}
\newcommand{\Ncal}{{\mathcal N}}

\newcommand{\loc}{{\text{\rm loc}}}

\newcommand{\Ecal}{{\mathcal E}}

\newcommand{\Ucal}{{\mathcal U}}

\newcommand{\rd}{{\text{\rm d}}}
\newcommand{\rw}{{\text{\rm w}}}

\newcommand{\rb}{{\text{\rm b}}}

\newcommand{\supp}{\operatorname*{supp}}

\newcommand{\Vinner}[1]{(\!({#1})\!)}
\newcommand{\Lim}{\operatorname*{\textsc{Lim}}_{T\rightarrow \infty}}

\begin{document}
\numberwithin{equation}{section}


\title[Convergence of time averages for the 3D NSE]{
  Convergence of time averages of weak solutions of the three-dimensional Navier-Stokes equations}

\author[C. Foias]{Ciprian Foias${}^1$}
\author[R. Rosa]{Ricardo M. S. Rosa${}^2$}
\author[R. Temam]{Roger M. Temam${}^3$}

\address{${}^1$ Department of Mathematics, Texas A\&M University,
  College Station, TX 77843, USA.}
\address{${}^2$ Instituto de Matem\'atica, Universidade Federal do Rio de Janeiro,
  Caixa Postal 68530 Ilha do Fund\~ao, Rio de Janeiro, RJ 21945-970,
  Brazil.}
\address{${}^3$Department of Mathematics, Indiana University, Bloomington,
  IN 47405, USA}

\email[R. Rosa]{rrosa@ufrj.br}
\email[R. M. Temam]{temam@indiana.edu}

\date{November 4, 2014}


\maketitle

\vspace{\baselineskip}

{\footnotesize \noindent \textsc{Abstract.}
Using the concept of stationary statistical solution, which generalizes the notion of invariant measure, it is proved that, in a suitable sense, time averages of almost every Leray-Hopf weak solution of the three-dimensional incompressible Navier-Stokes equations converge as the averaging time goes to infinity. This system of equations is not known to be globally well-posed, and the above result answers a long-standing problem, extending to this system a classical result from ergodic theory. It is also shown that, from a measure-theoretic point of view, the stationary statistical solution obtained from a generalized limit of time averages is independent of the choice of the generalized limit. Finally, any Borel subset of the phase space with positive measure with respect to a stationary statistical solution is such that for almost all initial conditions in that Borel set and for at least one Leray-Hopf weak solution starting with that initial condition, the corresponding orbit is recurrent to that Borel subset and its mean sojourn time within that Borel subset is strictly positive.

}

\section{Introduction}
An important quantity in the study of the asymptotic behavior of an evolutionary system, and in particular of the three-dimensional Navier-Stokes equations and of turbulent phenomena, is the \textbf{time average}
\[ \frac{1}{T} \int_0^T \Phi(\bu(t))\;\rd t
\]
of a solution $\bu=\bu(t)$, $t\geq 0$, where $T>0$ and $\Phi=\Phi(\bu)$ is a bounded Borel-measurable real-valued function defined on the phase space of the system and representing some information on the state of the system.

A fundamental question about these time averages is whether they converge as $T\rightarrow \infty$. This question can be raised independently of any kind of ergodic assumption for the flow.  Even for non-ergodic flows, an answer lies in the Pointwise Birkhoff Ergodic Theorem, which says that the time averages converge for almost every initial condition, with respect to any invariant measure of the flow, when the system is globally well-posed. If the invariant measure is ergodic, then these limits are the same for any solution, otherwise they may differ, but still converge almost everywhere. 

However, this question is quite more delicate for the three-dimensional Navier-Stokes equations, and for any other system which shares the same pathology of lacking a known result of global well-posedness. Without a well-defined semigroup at hand, invariant measures do not make sense, and the classical Birkhoff Theorem cannot be applied in the usual way. For this reason, the issue of convergence of the time averages for the three-dimensional Navier-Stokes equations has been an open problem for many decades. 

Nevertherless, in the conventional theory of turbulence, backed up by strong experimental evidence in fully-developed turbulent flows, the convergence of time averages and its equivalence to ensemble averages and, in some cases, to local averages in space are usually taken for granted (see e.g. \cite{taylor1935, taylor1938, batchelor1953, hinze1975, moninyaglom1975, frisch1995, lesieur1997}).

For a rigorous mathematical foundation for the study of turbulent phenomena and motivated by the lack of a well-posedness result for the Navier-Stokes equations, the notion of \textbf{statistical solution} was introduced in the early 1970's by Foias and Prodi \cite{foias72,foias73} (see also \cite{foiasprodi76}; G. Prodi does not appear as a co-author in the earlier works \cite{foias72,foias73}, but it is acknowledge in these papers that he was a major proponent of the theory), and a variant of this notion was introduced a few years later by Vishik and Fursikov \cite{vishikfursikov78}. The aim was to model the evolution of the probability distribution of the state of the system without resorting to an underlying semigroup. The part of the statistical information that does not change in time is captured by the so-called \textbf{stationary statistical solution}, which is a single measure representing the stationary statistics of the flow. In a sense, this notion generalizes the concept of invariant measure to the case in which global existence of individual solutions holds but global uniqueness does not necessarily hold. 

It is then natural to ask whether a pointwise Birkhoff-type result for the time averages holds with respect to stationary statistical solutions. We finally address this problem here for a class of stationary statistical solutions that we introduced recently in \cite{frt2010a, frt2010b, frtssp1}, bridging the two notions of statistical solutions given earlier in \cite{foias72, foias73, foiasprodi76, vishikfursikov78}. 

More precisely, this is done in the context of the three-dimensional incompressible Navier-Stokes equations on a bounded smooth domain with no-slip boundary condition. Moreover, the volume force is assumed to be steady and belonging to the phase space $H$ of square-integrable divergence-free velocity fields satisfying the corresponding boundary conditions. Similar results can be obtained for other boundary conditions. We let $H_\rw$ denote the space $H$ endowed with its weak topology and let $\Ccal_\rw=\Ccal_\loc([0,\infty),H_\rw)$ be the space of continuous functions from $[0,\infty)$ into $H_\rw$, endowed with the topology of uniform convergence on compact subintervals of $[0,\infty)$ with values in $H_\rw$ (see Section \ref{prelim} for more on the mathematical framework). We consider stationary statistical solutions associated with what we called in \cite{frt2010b} an \textbf{invariant (space-time) Vishik-Fursikov measure}, which is a Borel probability measure on the space $\Ccal_\rw$ carried by the set of Leray-Hopf weak solutions and invariant for the semigroup of translations in time (see Definition \ref{defvfinvmeas}).  Our convergence results below are valid almost everywhere with respect to such a measure. This result is reflected in the phase space in connection with the so-called \textbf{Vishik-Fursikov stationary statistical solution}, which is the projected measure $\rho_0=\Pi_0\rho$ on $H_\rw$ of an invariant (space-time) Vishik-Fursikov measure. The operator $\Pi_0$ is the projection map given by $\Pi_0\bu= \bu(0)$, i.e. the projection, onto $H_\rw$, of the value of a function $\bu$ in $\Ccal_\rw$ at time $t=0$ (see \eqref{defproj} and Definition \ref{defvfinvmeas}). 

Our first and key result consists in extending the Birkhoff Theorem to this situation, thus showing that time averages converge almost everywhere in a suitable sense. This result is the following one:

\begin{thm}
  \label{timeaverageclassicallimit}
  Let $\rho$ be an invariant Vishik-Fursikov measure on $\Ccal_\rw$ and let $\Phi$ be a bounded Borel-measurable real-valued function defined on the space $H$. Then, for $\rho$-almost every Leray-Hopf weak solution $\bu$, the classical limit
\begin{equation}
  \label{timeavelimit}
  \lim_{T\rightarrow \infty} \frac{1}{T} \int_0^T \Phi(\bu(t))\;\rd t
\end{equation}
exists. In particular, if $\rho_0=\Pi_0\rho$ is the associated Vishik-Fursikov stationary statistical solution, then for $\rho_0$-almost every initial condition $\bu_0$ in $H$, there exists at least one Leray-Hopf weak solution $\bu$ with $\bu(0)=\bu_0$ for which the limit in \eqref{timeavelimit} exists.
\end{thm}

Note that, in Theorem \ref{timeaverageclassicallimit}, the set of full measure of weak solutions for which the time averages converge may depend on the chosen invariant measure and on $\Phi$ (see, however, Remark \ref{indepentPhirmk}). We should emphasize, though, that \emph{any} invariant measure and \emph{any} such bounded function $\Phi$ are allowed; it is just the set of initial conditions for which the limit exists that may vary with these choices. Moreover, since a Leray-Hopf weak solution is actually bounded in $H$, the result of Theorem \ref{timeaverageclassicallimit} can be shown to hold for functions $\Phi$ which are not bounded on $H$, provided they are bounded on bounded subsets of $H$ (just multiply a non-bounded $\Phi$ by a truncation function that vanishes for sufficiently large values in $H$). This allows, for instance, for a function $\Phi$ representing the total kinetic energy of the flow on the domain $\Omega$ (the total kinetic energy is the density of the fluid times half the square of the norm in $H$ of the velocity field). The limitation for further applications of the result to experimentally relevant functions depend solely on the regularity of the Leray-Hopf weak solutions, which is a major standing problem in the mathematical theory of the Navier-Stokes equations.

\emph{Theorem \ref{timeaverageclassicallimit} is a positive answer to the long-standing problem of almost everywhere convergence of time averages of the incompressible three-dimensional Navier-Stokes equations. This has been the hope, as well as the belief, of our long-time collaborator, Oscar Manley, who insisted that we should address this issue. We are thus happy to finally answer this question.}

On a different, but related, matter, a particular stationary statistical solution can be constructed from the time averages of a weak solution chosen beforehand. This seems to have been constructed first (in the context of the Navier-Stokes equations) in \cite{foiastemam1975} through the use of convergent subsequences. The construction was improved in \cite{bcfm1995} (in the two-dimensional case) using the notion of generalized limit, which is a positive linear functional extending the classical limit to the space of essentially bounded functions (see Section \ref{prelim}). Given a weak solution $\bu$ and a generalized limit $\Lim$, there exists a stationary statistical solution $\mu_\bu$ on $H$ such that
\begin{equation}
  \label{defmusubu}
  \Lim \frac{1}{T} \int_0^T \Phi(\bu(t))\;\rd t = \int_H \Phi(\bv)\;\rd\mu_\bu(\bv),
\end{equation}
for any bounded Borel-measurable real-valued function $\Phi$ defined on $H$ (see Proposition \ref{timeavemeasure}). This construction is akin to the Bogoliubov-Krylov construction of invariant measures \cite{kb1937}. See \cite{wang2009,lukaszewickrealrobinson2011,chekrounglattholtz2012}, where this approach with generalized limits has been used to yield invariant measures to a large class of dissipative systems. See also  \cite{dapratozabczyk1996} for the corresponding construction in the context of stochastic equations.

Besides the issue of the generalized convergence being replaced by a classical convergence, addressed by Theorem \ref{timeaverageclassicallimit}, it was not known whether or not two different choices of generalized limits could yield two different stationary statistical solutions for the same weak solution. Using Theorem \ref{timeaverageclassicallimit}, we show that this cannot happen in general and that, at least almost everywhere with respect to any invariant Vishik-Fursikov measure, the choice of the generalized limit operator $\Lim$ is irrelevant from a measure-theoretic point of view. This result is stated in the following way:

\begin{thm}
  \label{timeaveragemeasureindepentofgenlim}
  Let $\rho$ be an invariant Vishik-Fursikov measure and let $\Lim$ and $\Lim'$ be two generalized limits. Then, for $\rho$-almost every Leray-Hopf weak solution $\bu$, the time-average stationary statistical solutions $\mu_\bu$ and $\mu_\bu'$ associated with $\bu$ and with $\Lim$ and $\Lim'$ according to \eqref{defmusubu}, respectively, are identical, i.e. $\mu_\bu=\mu_\bu'$ on $H$.
\end{thm}

Theorem \ref{timeaveragemeasureindepentofgenlim} is not an immediate consequence of Theorem \ref{timeaverageclassicallimit} since the set of measure zero in Theorem \ref{timeaverageclassicallimit} depends on the choice of $\Phi$; but using a topological separability argument this dependency can be overcome and Theorem \ref{timeaverageclassicallimit} can be used to prove Theorem \ref{timeaveragemeasureindepentofgenlim}. 

The time averages and Theorem \ref{timeaverageclassicallimit} can also be used to yield a strong form of recurrence. In fact, by taking $\Phi=\chi_E$ as the characteristic function of a Borel subset $E$ of $H$, the limit of the time averages gives the \textbf{mean sojourn time} of a weak solution $\bu$ in the set $E$:
\begin{equation}
  \label{defsoj}
  \text{Soj}_E(\bu)  = \lim_{T\rightarrow \infty} \frac{\left|\{ t\in [0,T]; \; \bu(t) \in E\}\right|}{T} =  \lim_{T\rightarrow \infty} \frac{1}{T} \int_0^T \chi_E(\bu(t))\;\rd t,
\end{equation}
where $|\cdot|$ denotes the Borel measure of a Borel subset of the real line. The application of Theorem \ref{timeaverageclassicallimit} in this case yields that if $E$ has positive $\rho_0$ measure (otherwise there is nothing of interest to prove), where $\rho_0=\Pi_0\rho$, then for $\rho$-almost every weak solution with $\bu(0)\in E$, the mean sojourn time $\text{Soj}_E(\bu)$ is strictly positive, meaning that the orbit is recurrent to $E$ and it returns to $E$ so often as to render the relative time that it spends within $E$ asymptotically positive:

\begin{thm}
  \label{sojournrecurrencethm}
  Let $\rho$ be an invariant Vishik-Fursikov measure and let $\rho_0=\Pi_0\rho$ be the associated Vishik-Fursikov stationary statistical solution. Let $E$ be a Borel subset of $H$ and assume $\rho_0(E)>0$. Then, for $\rho$-almost every Leray-Hopf weak solution $\bu$ with $\bu(0)\in E$, the sojourn time in $E$ of the orbit $\bu=\bu(t)$ is well-defined and is positive, i.e. $\text{Soj}_E(\bu)>0$. In particular, for $\rho_0$-almost every $\bu_0\in E$, there exists at least one Leray-Hopf weak solution $\bu$ on $[0,\infty)$ with $\bu(0)=\bu_0$ and such that the sojourn time in $E$ of the orbit $\bu=\bu(t)$ is positive, i.e. $\text{Soj}_E(\bu)>0$.
\end{thm}

Section \ref{prelim} is devoted to describing the precise framework in which the results are obtained, and the subsequent sections are devoted to proving each of the above theorems.

\section{Preliminares} 
\label{prelim}

We consider the incompressible Navier-Stokes equations with no-slip boundary conditions and a steady volume force $\bbf$, on a bounded spatial domain $\Omega\subset \RR^3$ with smooth boundary $\partial \Omega$. Denoting the kinematic viscosity by $\nu>0$, the velocity field by $\bu=\bu(t,\bx)$, and the kinematic pressure by $p=p(t,\bx)$, where $\bx$ and $t$ are the independent spatial and time variables, we write the equations as
\begin{equation}
  \label{nse}
  \begin{cases} \displaystyle
   \frac{\partial \bu}{\partial t} - \nu \Delta \bu
       + (\bu\cdot\bnabla)\bu + \bnabla p = \bbf, & \text{ on } (0,\infty)\times \Omega\\
      \bnabla \cdot \bu = 0, & \text{ on } (0,\infty)\times \Omega, \\ \bu = 0, &  \text{ on } (0,\infty)\times \partial \Omega.
   \end{cases}
\end{equation}
We assume that $\bbf\in H$, where $H$ is the closure, in $L^2(\Omega)^3$, of the space of infinitely-differentiable vector fields $\bv$ with compact support in $\Omega$ and satisfying the divergence-free condition $\bnabla \cdot \bv = 0$. The corresponding closure of this space in $H_0^1(\Omega)^3$ is also considered and is denoted by $V$. 

The space $H$ is endowed with the inner product inherited from $L^2(\Omega)^3$; the space $V$ is endowed with the inner product inherited from $H_0^1(\Omega)^3$; and we let $H_\rw$ denote the space $H$ endowed with its weak topology. We identify $H$ with its dual and consider the dual $V'$ of $V$, so that $V\subset H \subset V'$, with dense, compact, and continuous injections. 

We also consider the space $\Ccal_\rw=\Ccal_\loc([0,\infty),H_\rw)$ of functions continuous from $[0,\infty)$ to $H_\rw$ and endowed with the topology of uniform convergence on compact subintervals of $[0,\infty)$ with values in $H_\rw$. We will use the notation $\Ccal_\rw$ in some places for convenience.

A Leray-Hopf weak solution of \eqref{nse} on the interval of interest $[0,\infty)$ is a distribution function $\bu$ belonging to the space $L^\infty_\loc(0,\infty;H)\cap L^2_\loc(0,\infty;V)$ satisfying the system \eqref{nse} in a weak sense involving test functions in $V$;  being strongly continuous in $H$ from the right at the initial time $t=0$; and also satisfying the following energy inequality in the distribution sense on $[0,\infty)$:
\begin{equation}
  \frac{1}{2}\frac{\rd}{\rd t} \|\bu\|_{L^2}^2 + \nu \| \nabla \bu\|_{L^2}^2 \leq \Vinner{\bbf,\bu}_{L^2}.
\end{equation}
A Leray-Hopf weak solution automatically belongs to the space $\Ccal_\rw$. See e.g. \cite{lady63,temam,constantinfoias,fmrt2001a} for the mathematical background and the classical existence results for the Navier-Stokes equations.

We define the projection operator, from the space-time function space to the phase space, which is clearly continuous:
\begin{equation}
  \label{defproj}
  \Pi_0 : \Ccal_\loc([0,\infty),H_\rw) \rightarrow H_\rw, \quad \Pi_0\bu = \bu(0), \;\forall \bu\in \Ccal_\loc([0,\infty),H_\rw).
\end{equation}

We also consider, for any $\tau\geq 0$, the translation, or time-translation, operator
\begin{equation}
  \label{sigmatdef}
  \begin{aligned}
    \sigma_\tau :\Ccal_\loc([0,\infty),H_\rw) & \;\rightarrow \;\Ccal_\loc([0,\infty),H_\rw) \\
       \bu & \;\mapsto \;\sigma_\tau\bu \define (\sigma_\tau\bu)(t) = \bu(t+\tau), \qquad \forall t\geq 0.
  \end{aligned}
\end{equation}
In the context of the Navier-Stokes equations, a similar translation operator was used by Sell \cite{sell1996}.  The family $\{\sigma_\tau\}_{\tau\geq 0}$ is a continuous semigroup of linear operators on $\Ccal_\rw$. The family of time-translations also yields a jointly continuous operator:
\begin{equation}
  \label{sigmadef}
  \begin{aligned}
    \sigma :[0,\infty)\times \Ccal_\loc([0,\infty),H_\rw) & \;\rightarrow \;\Ccal_\loc([0,\infty),H_\rw) \\
       (\tau,\bu) & \;\mapsto \;\sigma_\tau\bu.
  \end{aligned}
\end{equation}

Let $\Ucal\subset \Ccal_\rw$ denote the set of Leray-Hopf weak solutions of the Navier-Stokes equations on the interval $[0,\infty)$:
\[ \Ucal = \{ \bu \in \Ccal_\loc([0,\infty),H_\rw); \bu \text{ is a Leray-Hopf weak solution on } [0,\infty) \}.
\]

We have proved in \cite{frtssp1} that $\Ucal$ is a Borel subset of $\Ccal_\rw$. We also remark that $\Ucal$ may not be invariant by the semigroup $\{\sigma_\tau\}_{\tau\geq 0}$. This is due to the assumption that a Leray-Hopf weak solution on $[0,\infty)$ is strongly continuous at the origin, without necessarily being strongly continuous at positive times. Hence, once we shift a Leray-Hopf weak solution $\bu$ to $\sigma_\tau \bu$, with $\tau>0$, the initial time of the shifted trajectory is the time $\tau$ of $\bu$, which is not necessarily a point of strong continuity of $\bu$. Nevertheless, $\sigma_\tau\bu$ satisfies all the other requirements, and we may say that $\sigma_\tau u$ is a Leray-Hopf weak solution when restricted to the open interval $(0,\infty)$.

We define an invariant Vishik-Fursikov measure and a Vishik-Fursikov stationary statistical solution as follows.
\begin{defs}
  \label{defvfinvmeas}
  An \textbf{invariant (space-time) Vishik-Fursikov measure} for the Navier-Stokes equations \eqref{nse} is a Borel probability measure $\rho$ on $\Ccal_\rw$ which is carried by $\Ucal$, i.e. $\rho(\Ucal)=1$, and is invariant by the translation semigroup $\{\sigma_\tau\}_{\tau \geq 0}$, i.e. $\rho(\sigma_\tau^{-1}\Ecal) = \rho(\Ecal)$, for any $\tau\geq 0$ and any Borel subset $\Ecal \subset \Ccal_\rw$. The projection $\rho_0 = \Pi_0\rho$ is a Borel probability measure on $H$ which is called a \textbf{Vishik-Fursikov stationary statistical solution}.
\end{defs}

The definition of invariant (space-time) Vishik-Fursikov measure given in \cite{frt2010b} is slightly different, but a regularity result for more general time-dependent Vishik-Fursikov measures proved in \cite[Theorem 4.1]{frtssp1} implies that the two definitions are in fact equivalent. The measure $\rho_0=\Pi_0\rho$ is a particular type of stationary statistical solution in the sense originally given in \cite{foias72} (see \cite{frt2010b} for more on this). Any invariant (space-time) Vishik-Fursikov measure and its corresponding Vishik-Fursikov stationary statistical solutions are carried by compact sets (see \cite{frt2010b}) and, hence, they are regular in the sense of measure theory \cite{schwartz1973,aliprantisborder2006}, i.e. the measure of any Borel set can be approximated from below by compact subsets of the Borel set and from above by open supersets of the Borel set.

A particular type of invariant measure is obtained via generalized limit of time averages. We recall that a \textbf{generalized limit} is any positive bounded real-valued linear functional defined on the space $L^\infty(0,\infty)$ and which extends the classical limit (see \cite{dunfordschwartz}). The following result appears in \cite{frt2010b}. More details about invariant (space-time) Vishik-Fursikov measures and Vishik-Fursikov stationary statistical solutions will be presented in the forthcoming work [Foias, C., Rosa, R., Temam, R.: Properties of stationary statistical solutions of the three-dimensional Navier-Stokes equations], which is a continuation of \cite{frtssp1}, focusing on the stationary case.

\begin{prop}
  \label{timeavemeasure}
  Let $\bu$ be a weak solution and let $\Lim$ be a generalized limit. Then, for any bounded Borel-measurable real-valued function $\varphi$ defined on $\Ccal_\rw$, the time averages of $\varphi(\sigma_t\bu)$ on the time intervals $(0,T)$ are uniformly bounded in $T$, so that
\begin{equation}
  \label{varphistargenlim}
  \Lim \frac{1}{T} \int_0^T \varphi(\sigma_t\bu)\;\rd t 
\end{equation}  
is well-defined. Moreover, there exists an invariant (space-time) Vishik-Fursikov measure $\rho_\bu$ such that
\begin{equation}
  \label{eqvarphistarpre}
  \Lim \frac{1}{T} \int_0^T \varphi(\sigma_t\bu)\;\rd t = \int_{\Ccal_\rw} \varphi(\bv)\;\rd\rho_\bu(\bv).
\end{equation}
In particular, for any bounded Borel-measurable real-valued function $\Phi$ defined on $H$ and for $\varphi = \Phi\circ \Pi_0$, the Vishik-Fursikov stationary statistical solution given by $\mu_\bu = \Pi_0\rho_\bu$ satisfies
\begin{equation}
  \label{eqPhistarpre}
  \Lim \frac{1}{T} \int_0^T \Phi(\bu(t))\;\rd t = \int_H \Phi(\bxi)\;\rd\mu_\bu(\bxi).
\end{equation}
\end{prop}

Thanks to Proposition \ref{timeavemeasure}, given a bounded Borel-measurable real-valued function $\varphi$ defined on $\Ccal_\rw$, we define the function $\varphi^*:\Ucal\rightarrow \RR$ through the relation
\begin{equation}
  \label{eqvarphistar}
  \varphi^*(\bu) = \Lim \frac{1}{T} \int_0^T \varphi(\sigma_t\bu)\;\rd t = \int_{\Ccal_\rw} \varphi(\bv)\;\rd\rho_\bu(\bv),
\end{equation}
and, similarly, given a bounded Borel-measurable real-valued function $\Phi$ defined on $H$, we define the function $\Phi^*:\Ucal \rightarrow \RR$ through
\begin{equation}
  \label{eqPhistar}
  \Phi^*(\bu) = \Lim \frac{1}{T} \int_0^T \Phi(\bu(t))\;\rd t = \int_H \Phi(\bxi)\;\rd\mu_\bu(\bxi).
\end{equation}

\section{Outline of the proof of Theorem \ref{timeaverageclassicallimit}}
\label{sectimeaverageclassicallimit}

Given an invariant Vishik-Fursikov measure $\rho$, we observe that the translation semigroup $\{\sigma_\tau\}_{\tau \geq 0}$ induces a semigroup on the space of Borel-measurable functions $\varphi: H\rightarrow \RR$, taking $\varphi$ into the function denoted $\sigma_\tau\varphi$ and defined by
\[ (\sigma_\tau\varphi)(\bu) \define \varphi(\sigma_\tau \bu).
\]
Since $\sigma_\tau$ is continuous on $\Ccal_\rw$, the function $\sigma_\tau\varphi$ is also Borel measurable on $H$. In particular, this is valid for $\varphi\in L^1(\rho)$, which means that $\{\sigma_\tau\}_{\tau \geq 0}$ is a strongly-measurable semigroup in $L^1(\rho)$. Clearly, $\sigma_\tau$ takes $L^\infty(\rho)$ into itself, satisfying, in the operator norm, $\|\sigma_\tau\|_{L^\infty(\rho)} \leq 1$, $\forall \tau \geq 0$. Moreover, since $\rho$ is invariant for the induced semigroup, we also have that $\sigma_\tau$ takes $L^1(\rho)$ into itself, satisfying, in the operator norm, $\|\sigma_\tau\|_{L^1(\rho)} = 1$, $\forall \tau \geq 0$. Therefore, by interpolation, it follows that $\{\sigma_\tau\}_{\tau \geq 0}$ is a strongly-measurable semigroup of bounded operators in $L^p(\rho)$, for $1\leq p \leq \infty$.

In this case, the Pointwise Birkhoff Ergodic Theorem \cite[Theorem VIII.7.5]{dunfordschwartz} (see also \cite[Section 1.2]{krengel1985}) says that the classical limit
\begin{equation}
  \label{appliedbirkhoff}
  \lim_{T\rightarrow \infty} \frac{1}{T}\int_0^T \varphi(\sigma_t\bu)\;\rd t
\end{equation}
exists for every $\varphi\in L^p(\rho)$, $1\leq p < \infty$, and for $\rho$-almost every $\bu$ in $\Ccal_\rw$.

In particular, if $\Phi:H\rightarrow \RR$ is a bounded Borel-measurable function from $H$ into $\RR$, then $\Phi$ is also a Borel-measurable function on $H_\rw$ since the two Borel $\sigma$-algebras of $H$ and $H_\rw$ are the same (which is due to the fact that they are both generated by the closed balls in $H$, which are closed for both topologies; see \cite{fmrt2001a, frtssp1}). Then, since $\Pi_0$ is continuous, the composition $\varphi = \Phi\circ \Pi_0$ is a bounded Borel-measurable function on $\Ccal_\rw$ and it belongs to all $L^p(\rho)$. Notice also that, for $\bu\in \Ccal_\rw$, we have
\[ \varphi(\sigma_t\bu) = \Phi(\Pi_0\sigma_t\bu) = \Phi(\bu(t)), \quad \forall t\geq 0.
\]
Then, \eqref{appliedbirkhoff} applies and we have that the classical limit \eqref{timeavelimit} exists for every bounded Borel-measurable function $\Phi: H \rightarrow \RR$ and for $\rho$-almost every $\bu$ in $\Ccal_\rw$. Since $\rho$ is carried by the set of Leray-Hopf weak solutions, this means in particular that the classical limit exists for $\rho$-almost every Leray-Hopf weak solution, which is the case of interest. This completes the first part of the proof.

Now we let $\Ncal$ be a set of $\rho$-null measure for which the limit \eqref{appliedbirkhoff} exists for any $\bu\in \Ccal_\rw\setminus \Ncal$. Since $\rho$ is carried by $\Ucal$, we may assume, without loss of generality, that $\Ccal_\rw\setminus \Ncal\subset \Ucal$. Set $N = H_\rw \setminus (\Pi_0(\Ccal_\rw \setminus \Ncal))$. Clearly,
  \[  \Pi_0^{-1} N = (\Pi_0^{-1}H_\rw) \setminus (\Pi_0^{-1}\Pi_0(\Ccal_\rw\setminus \Ncal)) \subset \Ccal_\rw \setminus (\Ccal_\rw\setminus \Ncal) = \Ncal.
  \]
  Thus, 
  \[ \rho_0(N) = \rho(\Pi_0^{-1} N) \leq \rho(\Ncal) = 0.
  \]
   Moreover, $H_\rw\setminus N = \Pi_0(\Ccal_\rw\setminus \Ncal)$, so that for any $\bu_0\in H_\rw\setminus N$, there exists $\bu\in \Ccal_\rw\setminus \Ncal$ such that $\bu(0)=\bu_0$. Since $N$ is of $\rho_0$-null measure as we have just proved and since \eqref{appliedbirkhoff} exists for any $\bu$ in $\Ccal_\rw\setminus \Ncal$, this means that for $\rho_0$-almost-every initial condition $\bu_0$, there exists at least one Leray-Hopf weak solution $\bu\in \Ccal_\rw\setminus \Ncal\subset \Ucal$ with $\bu(0)=\bu_0$ and for which the limit \eqref{appliedbirkhoff} exists. In particular, so does the limit \eqref{timeavelimit}, which completes the proof of the theorem.

\begin{rmk}
  \label{rhostarmeasurable}
  Given a Vishik-Fursikov measure $\rho$ and a function $\varphi \in L^1(\rho)$, then for each $T$, the map 
  \[ \bu \mapsto \frac{1}{T} \int_0^T \varphi(\sigma_t\bu) \;\rd t
  \]
  can be regarded as the composition of three Borel maps (the average on $(0,T)$, composed with $\varphi$, composed with the translations in time); hence it is a Borel map from $\Ccal_\rw$ into $\RR$. Thus, the function $\varphi^*$ defined by \eqref{eqvarphistar} is $\rho$-almost-everywhere the limit of Borel functions, and, therefore, $\varphi^*$ is measurable with respect to the Lebesgue completion of $\rho$. If $\varphi$ is a bounded Borel real-valued function defined on $\Ccal_\rw$, then it is integrable with respect to any probability measure on $\Ccal_\rw$, so that $\varphi^*$ is measurable with respect to the Lebesgue completion of any Vishik-Fursikov measure.
\end{rmk}

\begin{rmk}
  Using \eqref{appliedbirkhoff}, the Lebesgue Dominated Convergence Theorem, and the invariance of the measure $\rho$, we obtain the following identity (which is classical in the context of Ergodic Theory) for the function $\varphi^*$ given by \eqref{eqvarphistar}:
  \begin{equation}
    \label{birkhoffmeanid}
    \int_{\Ccal_\rw} \varphi^*(\bu) \;\rd\rho(\bu) = \int_{\Ccal_\rw} \varphi(\bu) \;\rd\rho(\bu),
  \end{equation}
  for any bounded Borel-measurable real-valued function $\varphi$ defined on $\Ccal_\rw$. The fact that $\varphi^*$ is $\rho$-measurable has been discussed in Remark \ref{rhostarmeasurable}. Since $\varphi^*$ is associated with the information obtained from the limit of time averages and the right-hand-side of \eqref{birkhoffmeanid} is the expected value of $\varphi$ with respect to an arbitrary Vishik-Fursikov invariant measure $\rho$, the relation \eqref{birkhoffmeanid} can be used to extend to arbitrary Vishik-Fursikov invariant measures relations that have been obtained for limits of time averages of Leray-Hopf weak solutions. This is a version in three dimensions of the two-dimensional result proved in \cite{fjmr2002} (see also \cite{bcfm1995}).
\end{rmk}

\begin{rmk}
In classical ergodic theory, it follows from the Poincar\'e Recurrence Theorem that the support of invariant measures is made of points which are nonwandering (i.e. any neighborhood of a point in the support contains a point whose orbit returns to that neighborhood infinitely often; see \cite{krengel1985}). A similar result is true in our case, as we will see in Remark \ref{nonwanderingrmk}. It may thus seem that the $\rho$-almost everywhere and $\rho_0$-almost everywhere results of Theorem \ref{timeaverageclassicallimit} fail to apply to a number of weak solutions that may contain some transient dynamics. In this sense, it is important to note that if $\bu$ and $\bv$ are two weak solutions with the same asymptotics in the weak topology, i.e. such that 
\begin{equation}
  \label{tracking}
  \bu(t) -\bv(t) \rightarrow 0, \text{ in } H_\rw, \quad \text{as } t\rightarrow \infty, 
\end{equation}
and if $\Phi$ is also weakly continuous on $H$, then
\[ \lim_{T\rightarrow \infty} \frac{1}{T} \int_0^T \left(\Phi(\bu(t))-\Phi(\bv(t))\right)\;\rd t = 0.
\]
Therefore, at least for this kind of functions $\Phi$, if the limit of the time averages exists for one of the weak solutions, then it also exists for the other, and the limits are the same. This extends the result of existence of the limit of the time averages to the collection of Leray-Hopf weak solutions which are asymptotically tracked (in the sense of \eqref{tracking}) by the globally defined (in the past and in the future), nonwandering Leray-Hopf weak solutions belonging to the support of the invariant measure, and, in particular, to the weak global attractor of the three-dimensional Navier-Stokes equations.
\end{rmk}

\begin{rmk}
  \label{classesofstatsolrmk}
  We have mentioned several types of statistical solutions which are potentially different and some of them actually live in different function spaces. In fact, we distinguish two main notions of statistical solution. The first, which we call simply \textbf{statistical solution}, was introduced in \cite{foias72,foias73}: it is a family of Borel probability measures on the phase space $H$ of the system, parametrized by the time variable, and representing the evolution in time of the probability distribution of the spatial velocity field of the flow. The second notion, which we termed a \textbf{(space-time) Vishik-Fursikov measure}, was inspired by the work of Vishik and Fursikov \cite{vishikfursikov78}, and formulated in \cite{frt2010a,frtssp1}: it is a single measure on the space of trajectories $\Ccal_\rw$,  representing the whole space-time distribution of the velocity field. Every (space-time) Vishik-Fursikov measure can be projected to a family of measures on $H$ which is a statistical solution, but it is not known whether every statistical solution can be lifted to a (space-time) Vishik-Fursikov measure (see \cite{frt2010a,frtssp1}).
\end{rmk} 

\begin{rmk}
A (space-time) Vishik-Fursikov measure may display statistics that change with time, such as those associated with evolving, and, in particular, decaying turbulence, or that may not, which is the case of those associated with stationary turbulence. In the latter case, we have a particular type of (space-time) Vishik-Fursikov measure which is invariant by the time-translation operator and is termed an \textbf{invariant (space-time) Vishik-Fursikov measure}. A particular type of invariant (space-time) Vishik-Fursikov measure is obtained via generalized limit of time averages of Leray-Hopf weak solutions (as in \eqref{eqPhistarpre}). 
\end{rmk}
  
\begin{rmk} 
Similarly, a particular type of statistical solution is that in which the family of measures does not vary with time, so that the statistics of the flow are statistically stationary. Such a statistical solution is called a \textbf{stationary statistical solution}. Moreover, as mentioned in Remark \ref{classesofstatsolrmk}, every (space-time) Vishik-Fursikov measure can be projected to a family of measures on $H$ which is a statistical solution. This yields a particular type of statistical solution which we termed a \textbf{Vishik-Fursikov statistical solution}. For instance, the projection of an invariant (space-time) Vishik-Fursikov measure yields a particular type of stationary statistical solution which we termed a \textbf{Vishik-Fursikov stationary statistical solution}. Finally, a particular type of Vishik-Fursikov stationary statistical solution is obtained via generalized limits of time averages of Leray-Hopf weak solutions (as in  \eqref{eqvarphistarpre}). 
\end{rmk}

\section{Outline of the proof of Theorem \ref{timeaveragemeasureindepentofgenlim}}

We want to establish that $\mu_\bu=\mu_\bu'$, for $\rho$-almost every $\bu$. For that purpose, using that the space $H$ is metrizable, it suffices to prove (see \cite[Theorem 15.1]{aliprantisborder2006}) that
\begin{equation}
  \label{equivbymeans}
  \int_H \Phi(\bv)\;\rd\mu_\bu(\bv) = \int_H \Phi(\bv) \;\rd\mu_\bu'(\bv), \quad \forall \Phi \in \Ccal_\rb(H).
\end{equation}
One difficulty in establishing \eqref{equivbymeans}, however, is that $\Ccal_\rb(H)$ is not separable. We use then that the Borel sets for the weak and for the strong topologies of $H$ are the same (as discussed just after \eqref{appliedbirkhoff}) so that $\mu_\bu$ and $\mu_\bu'$ can be viewed as Borel probability measures on $H_\rw$. Moreover, we use that $\mu_\bu$ and $\mu_\bu'$ are carried by the weakly compact set $K=B_H(R_0)$, for $R_0$ sufficiently large, which follows from the fact that the ball $B_H(R_0)$ contains the weak global attractor, provided $R_0 \geq |\bbf|_{L^2}/{\nu\lambda_1}$ (see \cite[Theorem IV.4.2]{fmrt2001a} or \cite{frt2010b}). The measures $\mu_\bu$ and $\mu_\bu'$ can then be viewed as Borel probability measures on $K_\rw$, and we then use that $\Ccal_\rb(K_\rw)=\Ccal(K_\rw)$ is separable, i.e. there exists a countable dense set $\Scal=\{\Phi_k\}_{k\in \NN}$ in $\Ccal(K_\rw)$. Therefore, $\mu_\bu=\mu_\bu'$ if, and only if,
\begin{equation}
  \label{equivbymeanscount}
  \int_H \Phi_k(\bv)\;\rd\mu_\bu(\bv) = \int_H \Phi_k(\bv) \;\rd\mu_\bu'(\bv), \quad \forall k\in \NN.
\end{equation}

For each $\Phi_k$, it follows from Theorem \ref{timeaverageclassicallimit} proved above that there exists a set $\Ecal_k$ in $\Ccal_\rw$ of zero $\rho$-measure such that the limit
\[ \lim_{T\rightarrow \infty} \frac{1}{T} \int_0^T \Phi_k(\bu(t))\;\rd t
\]
exists for every $\bu\in \Ccal_\rw\setminus \Ecal_k$. Take $\Ecal = \cup_k \Ecal_k$, so that $\Ecal$ still has $\rho$-measure zero and the classical limit exists for every $\bu\in \Ccal_\rw\setminus \Ecal$ and every $\Phi_k$, $k\in \NN$.

Since each $\Phi_k$ is also a bounded, strongly Borel-measurable function on $H$, we have, by the construction of $\mu_\bu$ and $\mu_\bu'$, that
\begin{multline*} 
  \int_H \Phi_k(\bv)\;\rd\mu_\bu(\bv) = \Lim\frac{1}{T} \int_0^T \Phi_k(\bu(t))\;\rd t = \lim_{T\rightarrow \infty} \frac{1}{T} \int_0^T \Phi_k(\bu(t))\;\rd t \\
    = \Lim\!\!{}' \frac{1}{T} \int_0^T \Phi_k(\bu(t))\;\rd t = \int_H \Phi_k(\bv)\;\rd\mu_\bu'(\bv),
\end{multline*}
for every $k\in \NN$. This proves \eqref{equivbymeanscount} and we deduce that $\mu_\bu=\mu_\bu'$, completing the proof.

\begin{rmk}
  \label{indepentPhirmk}
  The separability argument described above can be used to prove the existence of a single set of full measure for which the convergence of the time averages holds independently of the choice of $\Phi$ within the space $\Ccal(H_\rw)$. The idea starts by writing the space of Leray-Hopf weak solutions as a countable union $\Ucal = \cup_{n\in \NN} \Ucal \cap\Ccal_\loc([0,\infty),B_H(nR_0)_\rw)$, with $R_0$ as above, and considering countable collections $\Scal^n = \{\Phi_k^n\}_{k\in \NN}$ of functions in $\Ccal(B_H(R_0)_\rw)$ dense in this space. We extend each $\Phi_k^n$ to a bounded Borel measurable function $\bar\Phi_k^n$ by simply setting it to zero outside the ball $B_H(nR_0)_\rw$. Next, for each $\bar\Phi_k^n$, there exists a set $\Ecal_k^n$ of null $\rho$-measure in $\Ucal$ such that the time-average limit of $\Phi_k^n(\bu(\cdot)) = \bar\Phi_k^n(\bu(\cdot))$ exists. We set $\Ecal = \cup_{n,k\in \NN} \Ecal_k^n$, which is still of null $\rho$-measure. Then, for any $\bu\in \Ucal\setminus \Ecal$ and any $\Phi\in \Ccal(H_\rw)$, we use that $\bu\in \Ucal \cap\Ccal_\loc([0,\infty),B_H(nR_0)_\rw)$, for some $n\in \NN$, and that $\Phi$ restricted to $B_H(nR_0)_\rw$ is bounded, since this ball is (weakly) compact, to find a sequence $\Phi_{k_j}^n$ that converges, in $\Ccal(B_H(nR_0)_\rw)$, to the restriction of $\Phi$ to $B_H(nR_0)_\rw$. This allows us to show, by triangulation with $\Phi_k^n(\bu(\cdot))$, that the time averages of $\Phi(\bu(\cdot))$ are Cauchy and, hence, converge as $T\rightarrow \infty$.
\end{rmk}

\section{Outline of the proof of Theorem \ref{sojournrecurrencethm}}
\label{proofofsojournrecurrencethm}

Let $E$ be as in the statement of Theorem \ref{sojournrecurrencethm} and set $\Ecal = \Pi_0^{-1}E$, which is a Borel subset of $\Ccal_\rw$. Consider the set $\Scal_\Ecal$ in which the mean sojourn time exists, i.e.
\[ \Scal_\Ecal = \left\{ \bu \in \Ccal_\loc([0,\infty),H_\rw); \; \exists \lim_{T\rightarrow \infty} \frac{1}{T}\int_0^T \chi_{\Ecal}(\sigma_t \bu) \;\rd t = \text{Soj}_E(\bu) \right\}.
\]
Note that $\chi_\Ecal(\sigma_t\bu) = \chi_E(\bu(t))$, so that the limit in the definition above is in fact $\text{Soj}_E(\bu)$. We split $\Scal_\Ecal$ into the two subsets
\[ \Scal_\Ecal^0 = \left\{ \bu \in \Scal_\Ecal; \; \text{Soj}_E(\bu) = 0 \right\}, \quad \Scal_\Ecal^+ = \left\{ \bu \in \Scal_\Ecal; \; \text{Soj}_E(\bu) > 0 \right\}.
\]
The statement of Theorem \ref{sojournrecurrencethm} that the sojourn time in $E$ is positive for $\rho$-almost every Leray-Hopf weak solution $\bu$ with $\bu(0)\in E$ is precisely that $\rho(\Ecal\setminus \Scal_\Ecal^+)=0$. This is what we aim to prove. 

Using \eqref{birkhoffmeanid}, we write the measure of the set $\Ecal \setminus \Scal_\Ecal^+$ as
\begin{equation}
  \label{measureaux}
  \rho(\Ecal \setminus \Scal_\Ecal^+) =  \int_{\Ccal_\rw} \chi_{\Ecal \setminus \Scal_\Ecal^+}(\bu) \;\rd\rho(\bu) = \int_{\Ccal_\rw} \chi_{\Ecal \setminus \Scal_\Ecal^+}^*(\bu) \;\rd\rho(\bu),
\end{equation}
where $\chi_{\Ecal \setminus \Scal_\Ecal^+}^*$ is obtained from $\chi_{\Ecal \setminus \Scal_\Ecal^+}$ as in \eqref{eqvarphistar}.

Now we use the decomposition of $\Ccal_\rw$ into $\Scal_\Ecal^0 \cup \Scal_\Ecal^+ \cup (\Ccal_\rw \setminus \Scal_\Ecal)$, and the fact that $\Ccal_\rw \setminus \Scal_\Ecal$ is of null $\rho$-measure (which follows from Theorem \ref{timeaverageclassicallimit}) to write \eqref{measureaux} as
\begin{equation}
  \label{measureaux2}
  \rho(\Ecal \setminus \Scal_\Ecal^+) =  \int_{\Scal_\Ecal^0} \chi_{\Ecal \setminus \Scal_\Ecal^+}^*(\bu) \;\rd\rho(\bu) + \int_{\Scal_\Ecal^+} \chi_{\Ecal \setminus \Scal_\Ecal^+}^*(\bu) \;\rd\rho(\bu).
\end{equation}

It is immediate to check that $\Scal_\Ecal^+$ is invariant by the translation semigroup, i.e. 
\[ \sigma_\tau \Scal_\Ecal^+ \subset \Scal_\Ecal^+, \quad \forall \tau \geq 0.
\]
(And likewise for $\Scal_\Ecal^0$, but this will not be needed.) Hence, for every $\bu\in \Scal_\Ecal^+$, we have that $\sigma_\tau\bu \in \Scal_\Ecal^+$ for all $\tau \geq 0$, or, in other words, $\sigma_\tau\bu \notin \Ecal\setminus \Scal_\Ecal^+$. Thus, 
\[ \lim_{T\rightarrow \infty} \frac{1}{T} \int_0^T \chi_{\Ecal\setminus \Scal_\Ecal^+}(\sigma_t \bu) \;\rd t = 0, \quad \forall \bu \in \Scal_\Ecal^+.
\]
Since this time-average limit is $\rho$-almost everywhere equal to $\chi_{\Ecal \setminus \Scal_\Ecal^+}^*$, it follows that this function vanishes on $\Scal_\Ecal^+$, so that the second integral in \eqref{measureaux2} is zero.

Next, for $\rho$-almost every $\bu\in \Scal_\Ecal^0$, we have by \eqref{appliedbirkhoff} that
\[ \chi_{\Ecal\setminus \Scal_\Ecal^+}^*(\bu) = \lim_{T\rightarrow \infty} \frac{1}{T} \int_0^T \chi_{\Ecal\setminus \Scal_\Ecal^+}(\sigma_t \bu) \;\rd t \leq \lim_{T\rightarrow \infty} \frac{1}{T} \int_0^T \chi_\Ecal(\sigma_t \bu) \;\rd t = 0.
\]
Thus, $\chi_{\Ecal\setminus \Scal_\Ecal^+}^*$ vanishes on $\Scal_\Ecal^0$ as well, so that the first integral in \eqref{measureaux2} is also zero. This completes the proof that $\rho(\Ecal \setminus \Scal_\Ecal^+) = 0$.

Now we set $\Ncal = \Ecal\setminus (\Scal_\Ecal^+\cap \Ucal) \subset (\Ecal \setminus \Scal_\Ecal^+) \cup (\Ecal\setminus \Ucal)$. Since $\rho$ is carried by $\Ucal$ and we have proved that $\rho(\Ecal \setminus \Scal_\Ecal^+) = 0$, it follows that $\rho(\Ncal)=0$. Set $N=E \setminus (\Pi_0(\Ecal\setminus \Ncal))$. Then
  \[  \Pi_0^{-1} N = (\Pi_0^{-1}E) \setminus (\Pi_0^{-1}\Pi_0(\Ecal\setminus \Ncal)) \subset \Ecal \setminus (\Ecal\setminus \Ncal) = \Ncal,
  \]
  so that
  \[ \rho_0(N) = \rho(\Pi_0^{-1}N) \leq \rho(\Ncal) = 0.
  \]
  Moreover, since $\Pi_0\Ecal = \Pi_0\Pi_0^{-1} E = E$, we find that
  \[ E\setminus N = E \setminus (E \setminus (\Pi_0(\Ecal\setminus \Ncal))) = E \cap \Pi_0(\Ecal\setminus \Ncal) = \Pi_0\Ecal \cap \Pi_0 (\Ecal\setminus \Ncal) = \Pi_0(\Ecal\setminus \Ncal).
  \]
  Thus, for any $\bu_0\in E\setminus N$, there exists at least one $\bu\in \Ecal\setminus \Ncal\subset \Scal_\Ecal^+\cap \Ucal$ with $\bu(0)=\bu_0$. For such a $\bu$, we have that $\bu$ is a Leray-Hopf weak solution on $[0,\infty)$ starting at $\bu(0)=\bu_0$ and with $\text{Soj}_E(\bu)>0$. This completes the proof of Theorem \ref{sojournrecurrencethm}.

\begin{rmk}
  \label{nonwanderingrmk}
  Let $\bu_0$ be a point in the support $\supp\rho_0$ (the smallest closed subset of full measure) of a Vishik-Fursikov stationary statistical solution $\rho_0$. Since $\rho_0$ is a regular measure, it follows that if $O$ is a weakly open neighboorhood of $\bu_0$, then $\rho_0(O)>0$ (this is essentially proved in \cite[Theorem 12.14]{aliprantisborder2006}, although their definition of support is slightly different). Then, applying Theorem \ref{proofofsojournrecurrencethm} to the set $E=O$ we find that for $\rho_0$-almost every initial condition $\bv_0$ in $O$, there exists at least one Leray-Hopf weak solution $\bv$ on $[0,\infty)$ with $\bv(0)=\bv_0$ and for which the mean sojourn time $\text{Soj}_O(\bv)$ is strictly positive. In particular, this means that $\bu_0$ is nonwandering in the sense that for any weak neighborhood $O$ of $\bu_0$ there exists at least one initial condition $\bv_0$ in $O$ and one Leray-Hopf weak solution on $[0,\infty)$ starting at the point $\bv_0$ which returns to $O$ infinitely often. This shows that support $\supp\rho_0$ is made of points which are nonwandering in this sense.
\end{rmk}

\section*{Acknowledgments}

This work was partly supported by the National Science Foundation under the grants NSF-DMS-1206438 and NSF-DMS-1109784, by the Research Fund of Indiana University, and by the CNPq, Bras{\i}lia, Brazil, under the grant 303654/2013-9 and the cooperation project 490124/2009-7.

\end{document}